\documentclass[]{article}
\usepackage{graphicx}
\usepackage{caption}
\usepackage{subcaption}
\usepackage{algorithmic}
\usepackage{algorithm}
\usepackage{amsmath}
\usepackage{amssymb}
\usepackage{amsfonts}
\usepackage{amsthm}
\usepackage{fullpage}
\usepackage{url}
\usepackage{color} 
\usepackage{authblk} 
\usepackage{float}

\def\dt{{\textrm{d}t}}

\title{Nonstationary signal decomposition for dummies}
\author{Antonio Cicone   \\
DISIM\\
Universit\`a degli Studi dell'Aquila\\
via Vetoio n.1\\ 67100, L'Aquila\\
and\\
Istituto Nazionale di Alta Matematica\\
Citt\`a Universitaria\\
P.le Aldo Moro 5\\ 00185, Roma\\
ITALY\\
antonio.cicone@univaq.it}

\begin{document}

\maketitle

\begin{abstract}
How can I decompose a nonstationary signal? What are the advantages of using the most recent methods available in the literature versus using classical methods like (short time) Fourier transform or wavelet transform? This paper tries to address these and other questions providing the reader with a brief and self contained survey on what and how to tackle the decomposition of nonstationary signals.

\end{abstract}

\section{Introduction}

Nonstationary signals are ubiquitous in real life. We can consider, for instance, a stock market index, the ECG of a pregnant woman, or the terrestrial magnetic field measured by a magnetometer. For all these signals there are two kind of problems we may want to address:
\begin{description}
  \item[\textbf{Q1}] What are the active frequencies at each instant of time?
  \item[\textbf{Q2}] How to decompose such signals into simpler components?
\end{description}

We will see how these two questions, which may sound apparently unrelated, are indeed two sides of the same medal: understanding the nonstationary behavior of the signal.

Before addressing the previous two questions it is important to review what are the \emph{instantaneous frequency} and the \emph{time frequency representation} of a signal.

The idea of instantaneous frequency is natural. All of us know the basic sine function $y=A \sin(2 \pi\phi t)$ where $A$ is the amplitude, $1/\phi$ the period and $\phi$ the frequency. What if these quantities vary over time? We end up having what is called a amplitude modulated and frequency modulated (AM FM) signal $z~=~A(t) \sin(2 \pi)\phi(t)t)$. The intuition suggests to consider $\phi(t)$ as the instantaneous frequency of $z$ at the instant of time $t$. The question now is how to compute such quantity when a closed form formula of the signal is not known.

\begin{figure}[H]
        \centering\begin{subfigure}[b]{0.6\textwidth}
                \centering
                \includegraphics[width=\textwidth]{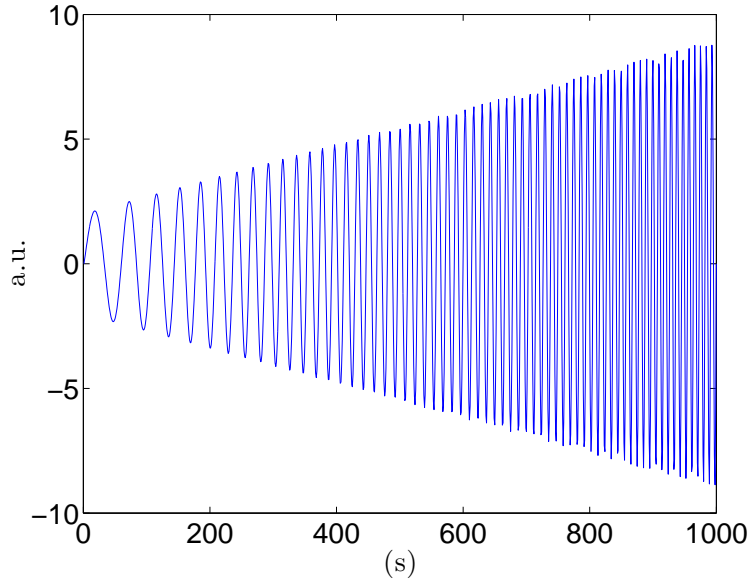}
        \end{subfigure}
        \caption{Example of a nonstationary signal with instantaneous frequency which
         increases over time. }\label{fig:Nonstationary_IMF}
\end{figure}

Ideally we would like to generalize the definition of stationary frequency as reciprocal of the period to the nonstationary case. However it becomes immediately clear that such generalization is, in most cases, not feasible. In Figure \ref{fig:Nonstationary_IMF} it is shown an example of a nonstationary AM FM signal with a frequency which increases over time. How to compute rigorously its instantaneous frequency? We cannot simply rely on the periods which cannot be computed exactly at each instant of time. One possible approach, well known and broadly diffused, is based on the evaluation of the Hilbert transform of the signal \cite[equation (3.4)]{huang1998empirical}. Another way, recently proposed in \cite[equation (34)]{cicone2014adaptive}, is based instead on the computation of the signal derivative. Both approaches have their own advantages as well as limitations. On the one hand the one based on the Hilbert transform, since it relies on integration, is inherently more stable, but at the same time the integration implies that not only local information are used for the instantaneous frequency computation. On the other hand the method  based on derivation is less stable, but at the same time it is completely local since only instantaneous information is used in its computation.

Once such information is made available we can represent it in a plot where the horizontal axis correspond to the time and the vertical one to the frequencies. This is what we call a time frequency representation.

The aforementioned methods allow to compute the instantaneous frequency of a signal whose instantaneous frequency is unique at each instant of time. 
What about signals which contain two or more instantaneous frequencies? In this case these two methods will not work.

This is where the questions Q1 and Q2 arise. In fact we can address the problem of studying a nonstationary signal in two ways. On the one hand we can try to compute directly its time frequency representation, addressing Q1. On the other hand we can tackle Q2 by first decomposing the signal into simple components each of which contains, ideally, a unique instantaneous frequency at each instant of time and then applying the previously mentioned techniques to compute the instantaneous frequencies component by component.

Historically researchers have first tackled Q1. The first technique used was the standard Fourier transform \cite{Bracewell1986} which is, however, inherently unable to capture any non stationarity in the signal. For this reason the so called short time Fourier transform was developed \cite{Cohen1995}. This last method allowed to produce a meaningful time frequency representations of nonstationary signals. The continuous wavelet transform \cite{Daubechies1992} allowed to improve the accuracy of the time frequency representation, in fact, instead on relying on an orthonormal basis of sines and cosines functions, it uses delation and translations of a mother wavelet.
To further sharpen the time frequency representation the so called synchrosqueezing trasform was proposed in \cite{Daubechies1996nonlinear,Daubechies2011Synchrosqueezed}.  Synchrosqueezing is a special kind of reassignment \cite{Flandrin1998,auger2013time}.

All the aforementioned methods are linear: they treat signals as linear combinations of elements in a preselected basis, sinusoidal or wavelet. There are also quadratic methods which are based instead on energy and power distributions, like the spectrogram, the scalogram, and the Wigner--Ville distribution. We refer the interested reader to the book \cite{Flandrin1998}.

It is important to remind at this point the so called Heisenberg uncertainty principle \cite{Cohen1995}. Based on this principle the accuracy that it can be achieved in producing a time frequency representation is limited. In particular the idea is that either we achive a good accuracy in assigning the frequency, but then the corresponding time horizon is not well identified, or vice versa \cite{flandrin2012uncertainty}.

In 1998, meanwhile new methods for a more accurate then ever time frequency representations were proposed, Huang and his research group at NASA devised the so called Empirical Mode Decomposition (EMD) algorithm \cite{huang1998empirical}. This method is the first technique ever developed able to address Q2 without any kind of a priori information and without making any a priori selection of a basis to be used. Furthermore, in some sense, it allows to bypass the Heisenberg uncertainty principle. In fact, decomposing directly a signal without any a priori information on its instantaneous frequencies allows to address Q2 avoiding the limitations induced by the Heisenberg uncertainty principle.

Few years later it became clear that this method is unstable. In particular its sensitivity to noise was unraveled. For this reason Huang and his group devised the so called Ensemble Empirical Mode Decomposition (EEMD) algorithm \cite{wu2009ensemble} which allows to overcome the sensitivity of the original EMD algorithm.

The publication of the EMD first and the EEMD after and their success inspired many other researchers to work on alternative methods for the decomposition of a signal into a few simple and meaningful components. All the alternative methods proposed so far are based on the minimization of some functional, like the sparse time--frequency representation algorithm \cite{hou2009variant,thomas2011adaptive}, and they require to make some assumption on the signal under study. Only one of them, called Iterative Filtering (IF) method \cite{lin2009iterative}, and its generalization, the Adaptive Local Iterative Filtering (ALIF) algorithm \cite{cicone2014adaptive}, are based on iterations like EMD and EEMD, and, therefore, no assumptions are required on the kind of signal we want to decompose.

The rest of this paper focuses on decomposition methods based on iterations for 1D signals. In particular we review the EEMD method, Section \ref{sec:EEMD},  
and the IF algorithm, Section \ref{sec:IF}. The paper ends with concluding remarks and an outlook to the main open problems in the field.

We point out here that both EEMD and IF have been generalized to 2D. See \cite{huang2009multi} and \cite{cicone2017multidimensional} for further details.


\section{The Ensemble Empirical Mode Decomposition algorithm}\label{sec:EEMD}

As we mentioned in the introduction, the goal of the EMD and EEMD methods is the decomposition of a signal into simple components, called Intrinsic Mode Functions (IMFs), each of which with a unique instantaneous frequency at every instant of time.

We start reviewing the EMD algorithm \cite{huang1998empirical}, whose pseudocode is given in Algorithm \ref{algo:EMD}.

\begin{algorithm}
\caption{\textbf{Empirical Mode Decomposition} IMF = EMD$(f)$}\label{algo:EMD}
\begin{algorithmic}
\STATE IMF = $\left\{\right\}$
\WHILE{the number of extrema of $f$ $\geq 2$}
      \STATE $s_1 = s$
      \WHILE{the stopping criterion is not satisfied}
                  \STATE  compute the moving average $M(s_{m}(x))$
                  \STATE  $s_{m+1}(x) = s_{m}(x) - M(s_m(x))$
                  \STATE  $m = m+1$
      \ENDWHILE
      \STATE IMF = IMF$\,\cup\,  \{ s_{m}\}$
      \STATE $s=s-s_{m}$
\ENDWHILE
\STATE IMF = IMF$\,\cup\,  \{ s\}$
\end{algorithmic}
\end{algorithm}

The key idea behind this method is what the authors called the sifting process: given a signal $s$ we capture its highest frequency oscillations by subtracting its moving average $M(s)$ from the signal itself. To do so we need to compute somehow the function $M(s)$. Huang and his research group proposed to compute first the upper and lower envelopes connecting the maxima and minima of the signal, respectively, by means of cubic splines. Then we compute the moving average as mean between these two curves point by point.

Since the derived moving average is only an approximation of the exact one the idea is to iterate the aforementioned calculation applying it to the new signal generated after the subtraction. In the end we expect the method to converge to an IMF or, using a stopping criterion, we discontinue the calculations when we are close in some sense to an IMF.

Using the approach just described, we compute the first IMF and we subtract it from the signal under study. The remainder can be treated as a new signal to which we apply again the sifting process. In the end we decompose the original signal into several IMFs and a remainder that cannot be decomposed anymore because it does not contain any oscillations.

Everything works fine except that the EMD process proved to be sensitive to small perturbations. In particular if we perturb a given signal with white noise, even if of small amplitude compared with the signal itself, the EMD may end up providing a completely different decomposition.

To address this issue Huang and his collaborators proposed in \cite{wu2009ensemble} to add to the given signal hundreds of white noise realizations. The final decomposition is then computed as the average of all these decompositions.

The EEMD Matlab implementation can be download from \url{http://rcada.ncu.edu.tw/}.

Besides the signal there are other two inputs that we need to pass to the EEMD Matlab function:
\begin{description}
  \item[Nstd] which represent the ratio between the standard deviation of the added noise and that of the signal under study.
  \item[NE] Number of noise realizations in the ensemble, id est the number of noise realizations to be added to the given signal.
\end{description}

The authors of the EEMD suggest in \cite{wu2009ensemble}  to set Nstd to 0.2. They also point out that if the data is dominated by high-frequency signals, the noise amplitude may be smaller, and when the data is dominated by low-frequency signals, the noise amplitude may be increased.

Furthermore they suggest to use ensembles of a few hundreds perturbations of a single signal. Based on our experience setting NE to one hundred is enough in most cases.

\subsection{Numerical examples}

Following what has been done in \cite{huang2008review} we apply the EEMD algorithm available at \url{http://rcada.ncu.edu.tw/} to two well known geophysical nonstationary signals: the Vostok temperature derived from ice core signal \cite{petit1997four,petit1999climate} and the length--of--day data (LOD) \cite{gross2000combinations}. For a detailed description about these two datasets and the meaning of their decompositions we refer the interested reader to \cite{huang2008review}.

\begin{figure}[H]
        \begin{subfigure}[b]{0.32\textwidth}
                \centering
                \includegraphics[width=\textwidth]{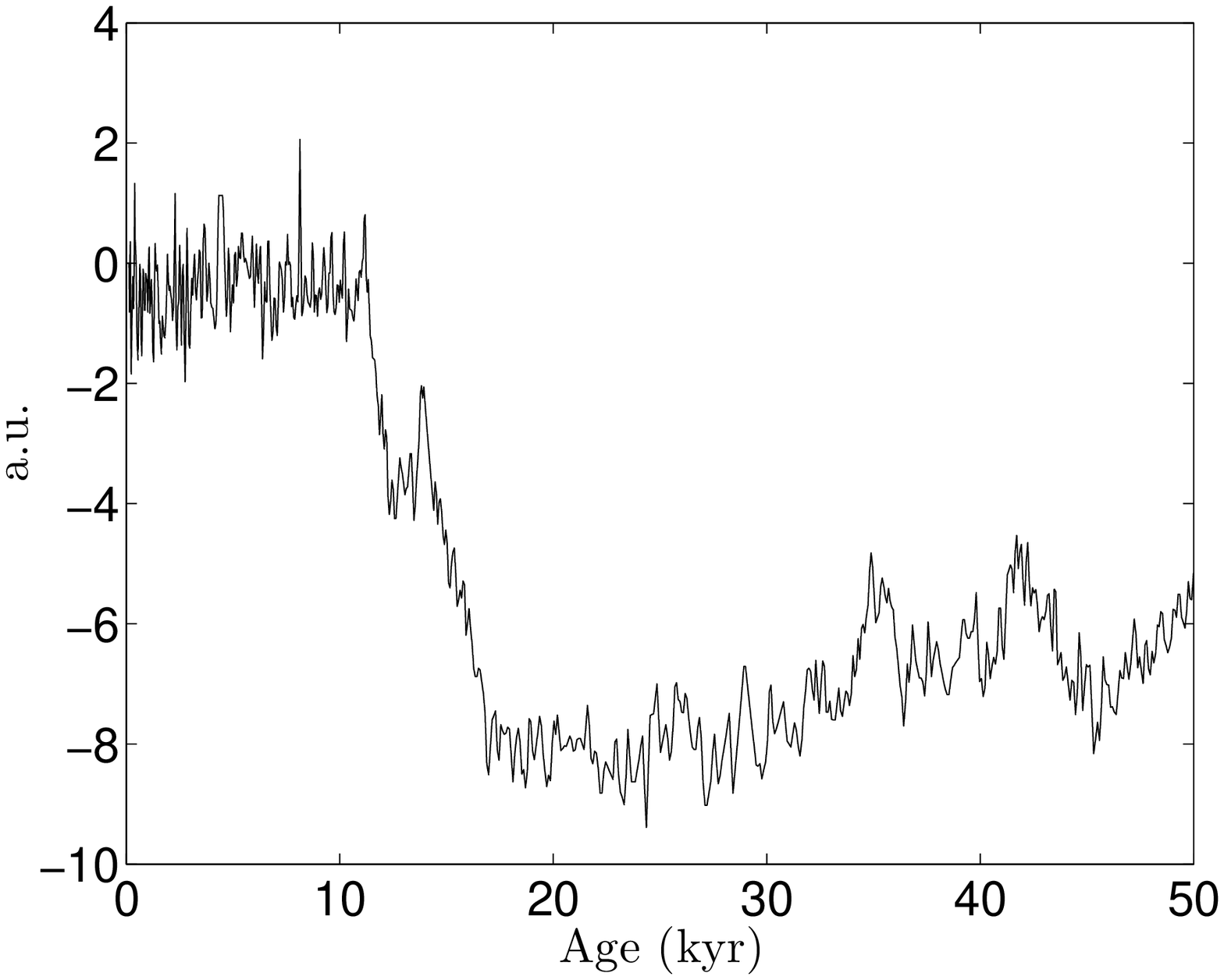}
                \caption{}
                \label{fig:Paleo}
        \end{subfigure}%
        ~ 
        \begin{subfigure}[b]{0.32\textwidth}
                \centering
                \includegraphics[width=\textwidth]{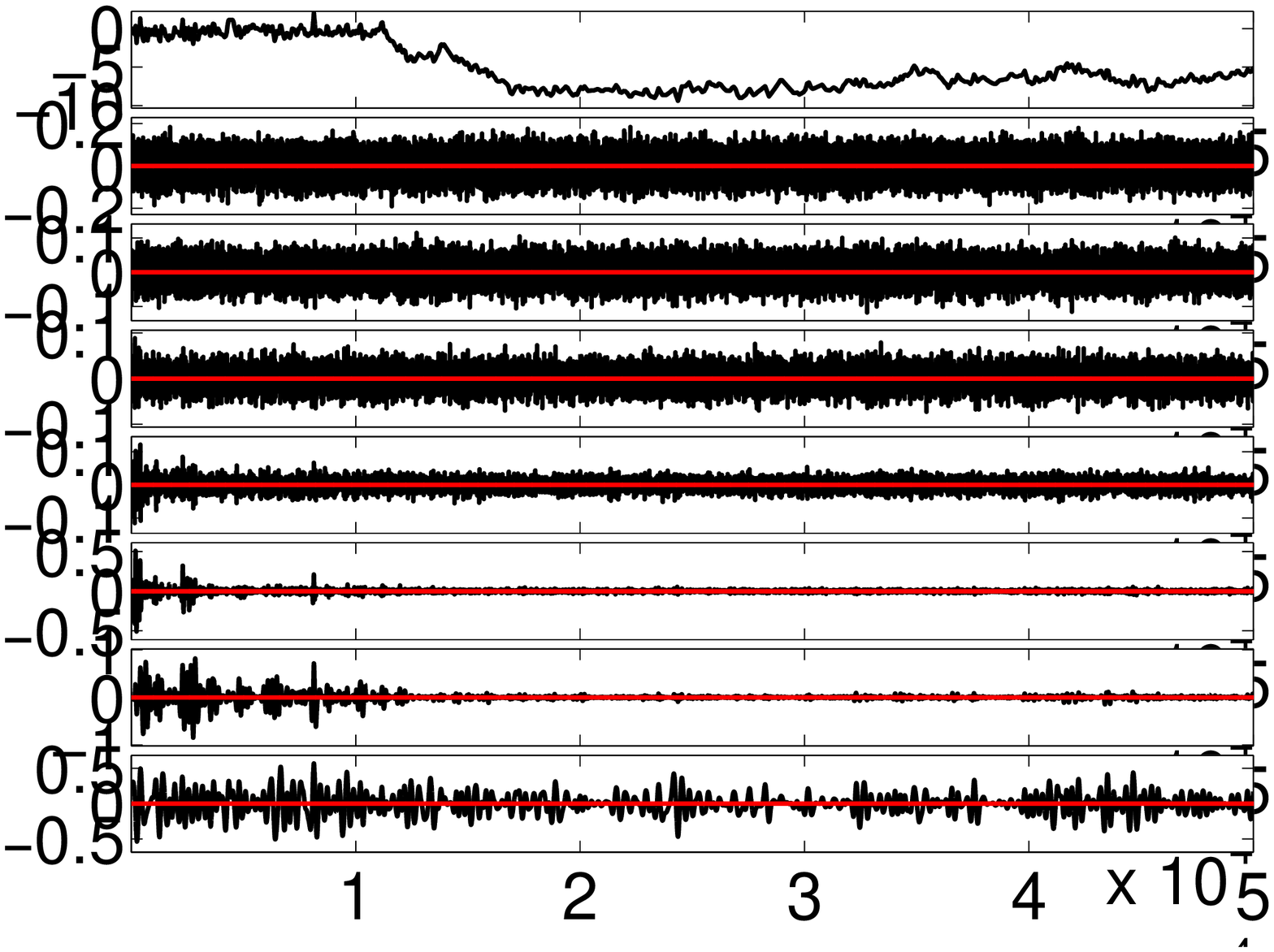}
                \caption{}
                \label{fig:Paleo_EEMD_Decomp_1}
        \end{subfigure}
        ~ 
        \begin{subfigure}[b]{0.32\textwidth}
                \centering
                \includegraphics[width=\textwidth]{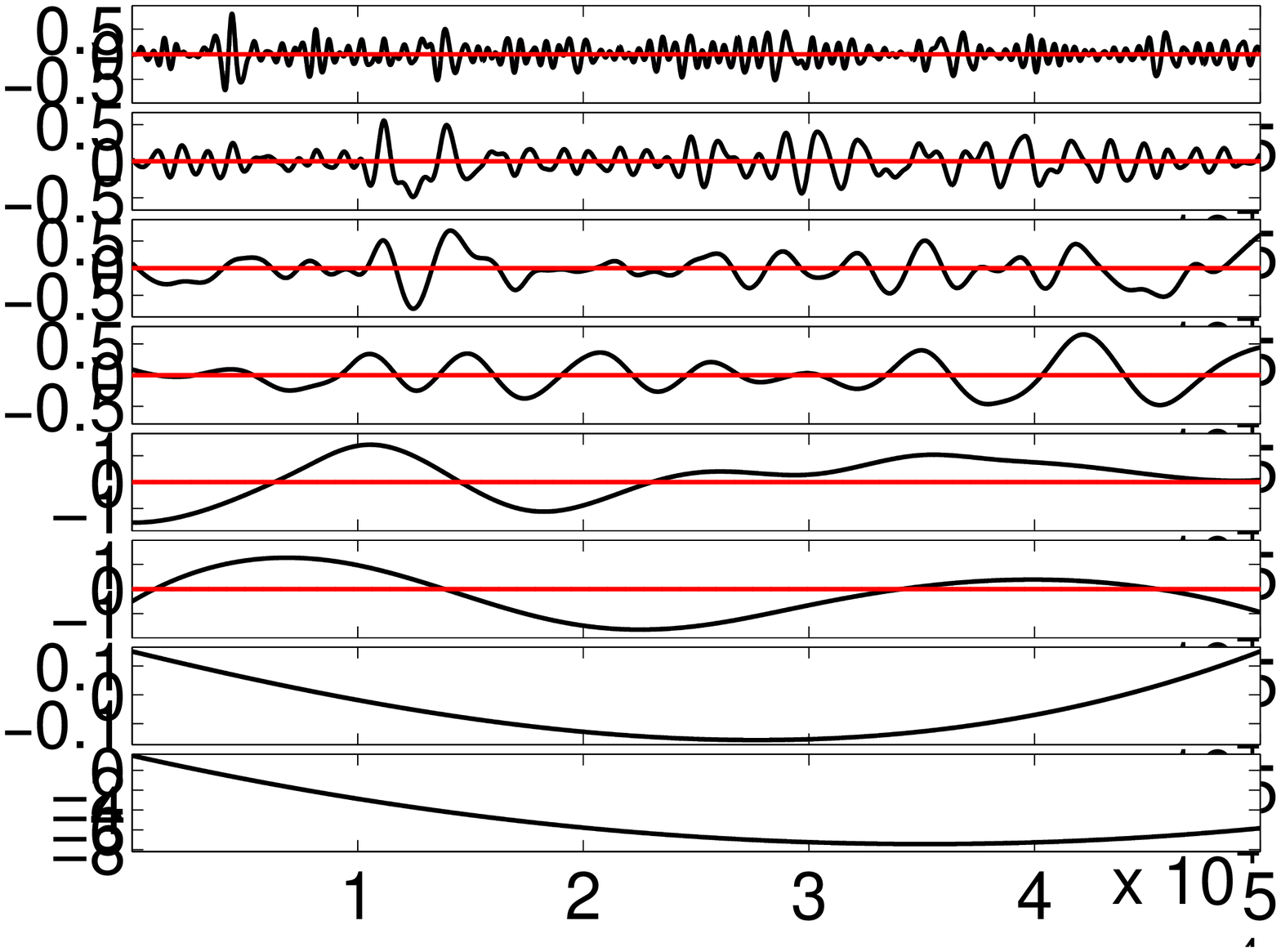}
                \caption{}
                \label{fig:Paleo_EEMD_Decomp_2}
        \end{subfigure}
        \caption{ (\subref{fig:Paleo}) Vostok temperature dataset of the last 50 thousands years. (\subref{fig:Paleo_EEMD_Decomp_1}) and (\subref{fig:Paleo_EEMD_Decomp_2}) EEMD decomposition. }\label{fig:Paleo_EEMD}
\end{figure}

We start with the Vostok temperature dataset. We focus, for simplicity, on the last 50 thousands years sample values, Figure \ref{fig:Paleo}. Applying the EEMD method with Nstd $=\ 0.2$, and NE $=\ 100$ we obtain the decomposition shown in Figures \ref{fig:Paleo_EEMD_Decomp_1} and \ref{fig:Paleo_EEMD_Decomp_2}.

\begin{figure}[H]
        \begin{subfigure}[b]{0.32\textwidth}
                \centering
                \includegraphics[width=\textwidth]{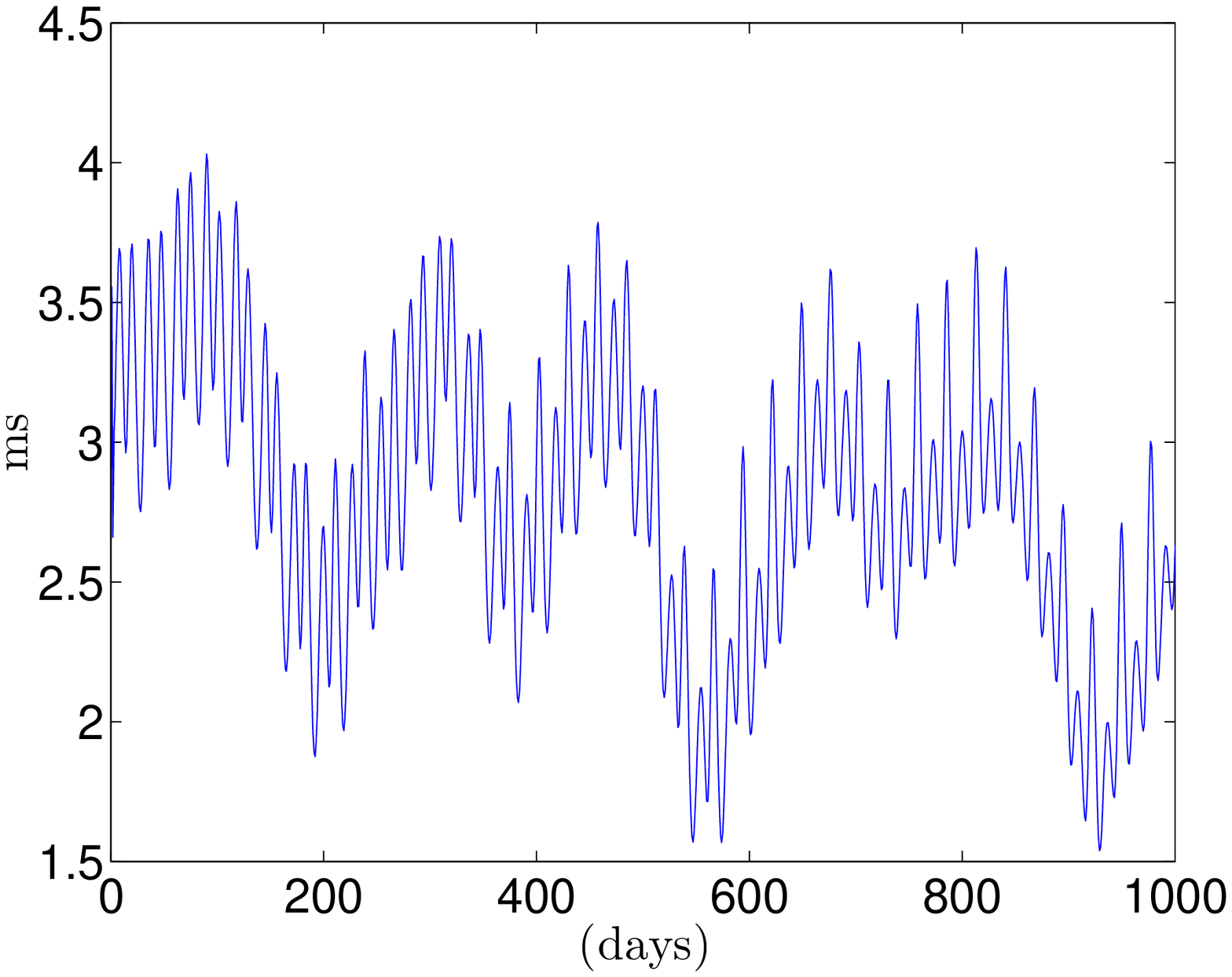}
                \caption{}
                \label{fig:LOD}
        \end{subfigure}%
        ~ 
        \begin{subfigure}[b]{0.32\textwidth}
                \centering
                \includegraphics[width=\textwidth]{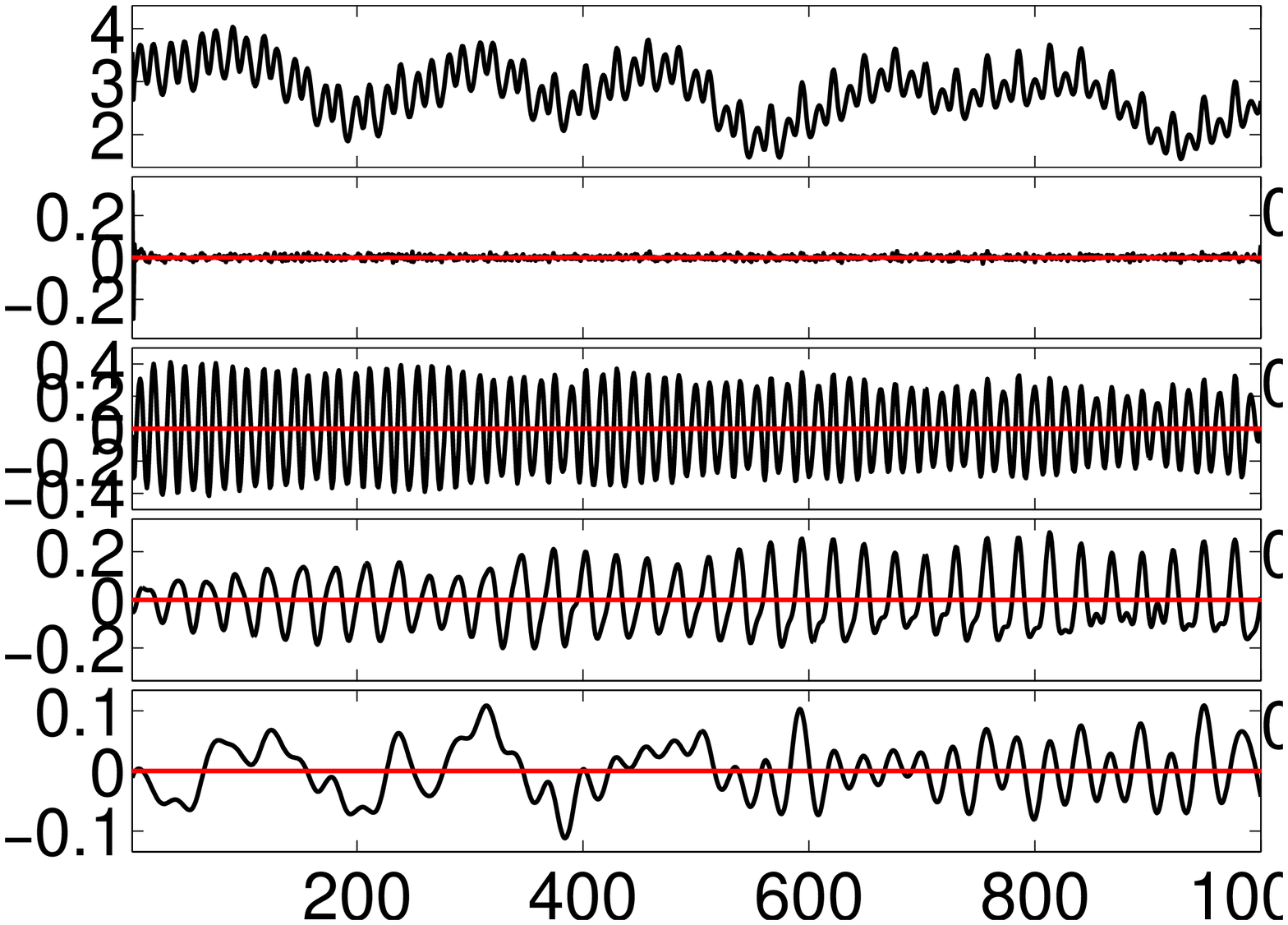}
                \caption{}
                \label{fig:LOD_EEMD_Decomp_1}
        \end{subfigure}
        ~ 
        \begin{subfigure}[b]{0.32\textwidth}
                \centering
                \includegraphics[width=\textwidth]{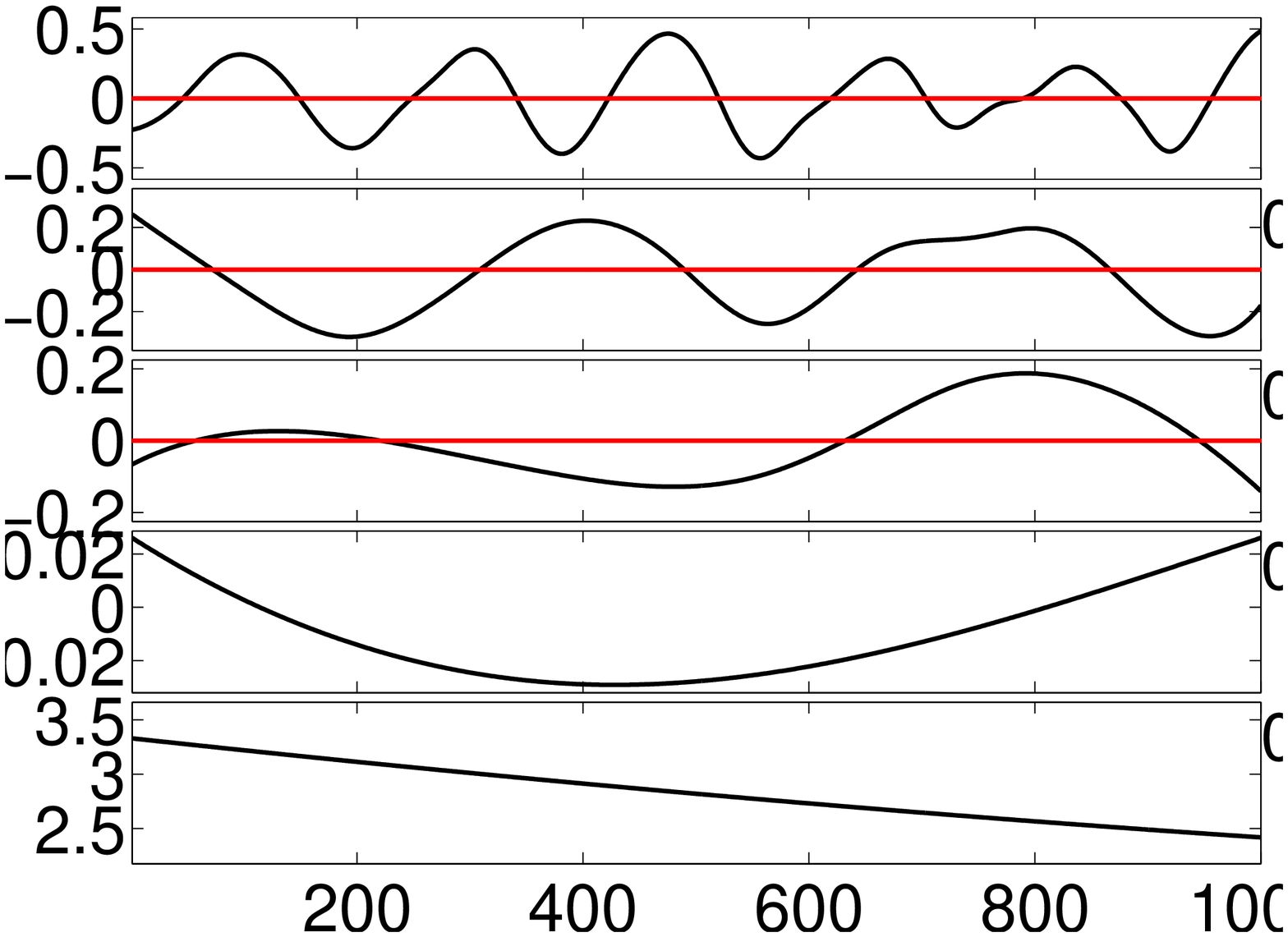}
                \caption{}
                \label{fig:LOD_EEMD_Decomp_2}
        \end{subfigure}
        \caption{ (\subref{fig:LOD}) length--of--day data from the beginning of 1983 to the end of 1986. (\subref{fig:LOD_EEMD_Decomp_1}) and (\subref{fig:LOD_EEMD_Decomp_2}) EEMD decomposition. }\label{fig:LOD_EEMD}
\end{figure}

For the LOD dataset, instead, we consider the data from the beginning of 1983 to the end of 1986. The signal is shown in Figure \ref{fig:LOD}.

If we apply the EEMD method with Nstd set to 0.2, and NE set to 100, we obtain the decomposition shown in Figures \ref{fig:LOD_EEMD_Decomp_1} and \ref{fig:LOD_EEMD_Decomp_2}.

\section{The Iterative Filtering method}\label{sec:IF}

Inspired by the EMD algorithm, Lin et al proposed in 2009 an alternative method called Iterative Filtering \cite{lin2009iterative}. The structure is the very same as in EMD with the only difference that now the moving average is computed as a local average of the values of the signal. This is achieved by integration of the signal itself weighted using an a priori chosen mask with non zero values concentrated on a finite interval.

How to choose such mask function? Based on the theoretical results \cite{cicone2014adaptive,cicone2017numerical,cicone2017spectral,cicone2017multidimensional} we know that to guarantee the convergence and stability of the algorithm, which in turn guarantee the physical meaningfulness of the decomposition, it is sufficient to consider a mask function which is generated as convolution with itself of a function fulfilling the following properties: compactly supported (it is zero outside a closed and finite set), nonnegative valued, with integral equal to 1, even (symmetric with respect to the vertical axis). In \cite{cicone2014adaptive} Fokker--Plank filters were proposed. They have the extra property of being infinitely smooth on the entire real line. The only drawback of such filters is that they are not known in an explicit form. However it is possible to compute them up to machine precision using numerical methods. One of such filters is available, together with a Matlab version of Iterative Filtering, online\footnote{\url{www.cicone.com}, GitHub and Mathworks}.

\begin{algorithm}
\caption{\textbf{Iterative Filtering} IMF = IF$(s)$}\label{algo:IF}
\begin{algorithmic}
\STATE IMF = $\left\{\right\}$
\WHILE{the number of extrema of $s$ $\geq N\geq 2$}
      \STATE $s_1 = s$
      \WHILE{the stopping criterion is not satisfied}
                  \STATE  compute the filter length $l_m$ for $s_{m}(x)$
                  \STATE  $s_{m+1}(x) = s_{m}(x) -\int_{-l_m}^{l_m} s_m(x+t)w_m(t)\dt$
                  \STATE  $m = m+1$
      \ENDWHILE
      \STATE IMF = IMF$\,\cup\,  \{ s_{m}\}$
      \STATE $s=s-s_{m}$
\ENDWHILE
\STATE IMF = IMF$\,\cup\,  \{ s\}$
\end{algorithmic}
\end{algorithm}

The pseudocode of the IF algorithm is given in Algorithm \ref{algo:IF} where $w_m(t)$ is the chosen mask function whose support is in $[-l_m,\ l_m]$, and $l_m$ is called \emph{mask length} which represents the half support length.

The current Matlab release of IF is \verb"IF_v6" and it requires as input the signal we want to decompose and, optionally, a variable, generated using an extra function called \verb"Settings_IF", which contains all the tuning parameters needed in the main algorithm. This implementation allows also the user to pass a third input, a vector containing a priori determined mask length values, forcing the method to skip their computation.

As outputs the algorithm returns a matrix containing as rows the IMFs, that we recall are the simple components in which the signal is decomposed, and a vector containing the values of mask length used to decompose each IMF.

What are the tuning parameters that can be set in IF by means of the Matlab function \verb"Settings_IF"? Let us review them one by one.

\begin{description}
  \item[delta] In the IF algorithm we need to use a stopping criterion to discontinue the calculation for each IMF, as pointed out in line four of the pseudocode. The one currently implemented, called \verb"delta", corresponds to the ratio between the norm 2 of the moving average curve and the norm 2 of the signal. Since the moving average curve converges towards a zero function we discontinue the calculations when the aforementioned ratio is small enough. Default value is set to 0.001.
  \item[ExtPoints] The algorithm iterates until the remainder has at most \verb"ExtPoints" number of extrema. This number corresponds to the value $N$ in line two of the pseudocode. Default value is equal to 3.
  \item[NIMFs] Maximal number of IMFs allowed in the decomposition, excluding the remainder. Default value is set to 1.

  \item[extensionType] The IF algorithm requires, as any other method in signal processing, to make an assumption on how the signal extends outside the boundaries. Three options are given: constant, periodical, and reflection. In the constant case the signal extends outside the boundaries with the last values achieved at the boundaries. Whereas periodical implies that the signal is assumed to repeat infinitely many times outside the boundaries. Finally reflection means that we extend the signal assuming symmetry with respect to the vertical lines passing through the boundary points.
   \item[MaxInner] Maximum number of iterations allowed for the computation of each IMF.  Default value is set to 200.
   \item[alpha] Parameter used for the mask length computation. In particular the algorithm measures the distances between subsequent extrema in the signal under study. It is up to the user to decide which value to be used in the decomposition. In fact, if \verb"alpha" is set to 0 the mask length is proportional to the minimum distance between two subsequent extrema. If it is set to 1 then it is proportional to the maximum distance. If we set it to \verb"ave" the mask length equals to, as suggested in \cite{lin2009iterative}, the roundoff value of the ratio $\frac{2*\textrm{Xi}*L}{N_e}$, where $L$ is the length of the signal, Xi is defined in the following, and $N_e$ represents the number of extrema in the signal. Finally if we set it to \verb"Almost_min" the mask length is selected equal to the round off value of $2*\textrm{Xi}*P_{30}$, where $P_{30}$ represents the 30 percentile of the vector containing all the distances between subsequent extrema of the signal. Default value set to \verb"ave".

       We point out that the smaller is the mask length the higher is the number of IMFs produced in the decomposition and the finer is the separation of two frequencies. Based on our experience the option \verb"Almost_min" allows to obtain a decomposition with a high number of components, but with IMFs which have in unique instantaneous frequencies all the time. Higher values in the parameter \verb"alpha" may lead, in some cases, to components which show an intermittency in the frequency.
  \item[Xi] This last parameter allows to tune the mask length. Depending on the chosen filter, in fact, we may need to extend the filter length more or less due to its shape. In particular if we select a mask shape with big values in the center, but which goes quickly to zero we may need to select \verb"Xi" bigger than 2. Whereas if the mask shape is almost constant everywhere on its support then \verb"Xi" may be selected closer to 1. Based on the numerical analysis conducted in \cite{cicone2017numerical} we know that enlarging the mask function support implies squeezing its fourier transform and in turn it allows to tune the sampled frequencies that enter in the extracted IMF. Suggested values range in the interval $[1.1,\ 3]$. Default value is equal to 1.6.
\end{description}

\subsection{Numerical examples}

We apply the \verb"IF_v6" to the Vostok temperature derived from ice core signal \cite{petit1997four,petit1999climate} and the length--of--day data (LOD) \cite{gross2000combinations}.

\begin{figure}[H]
        \begin{subfigure}[b]{0.5\textwidth}
                \centering
                \includegraphics[width=\textwidth]{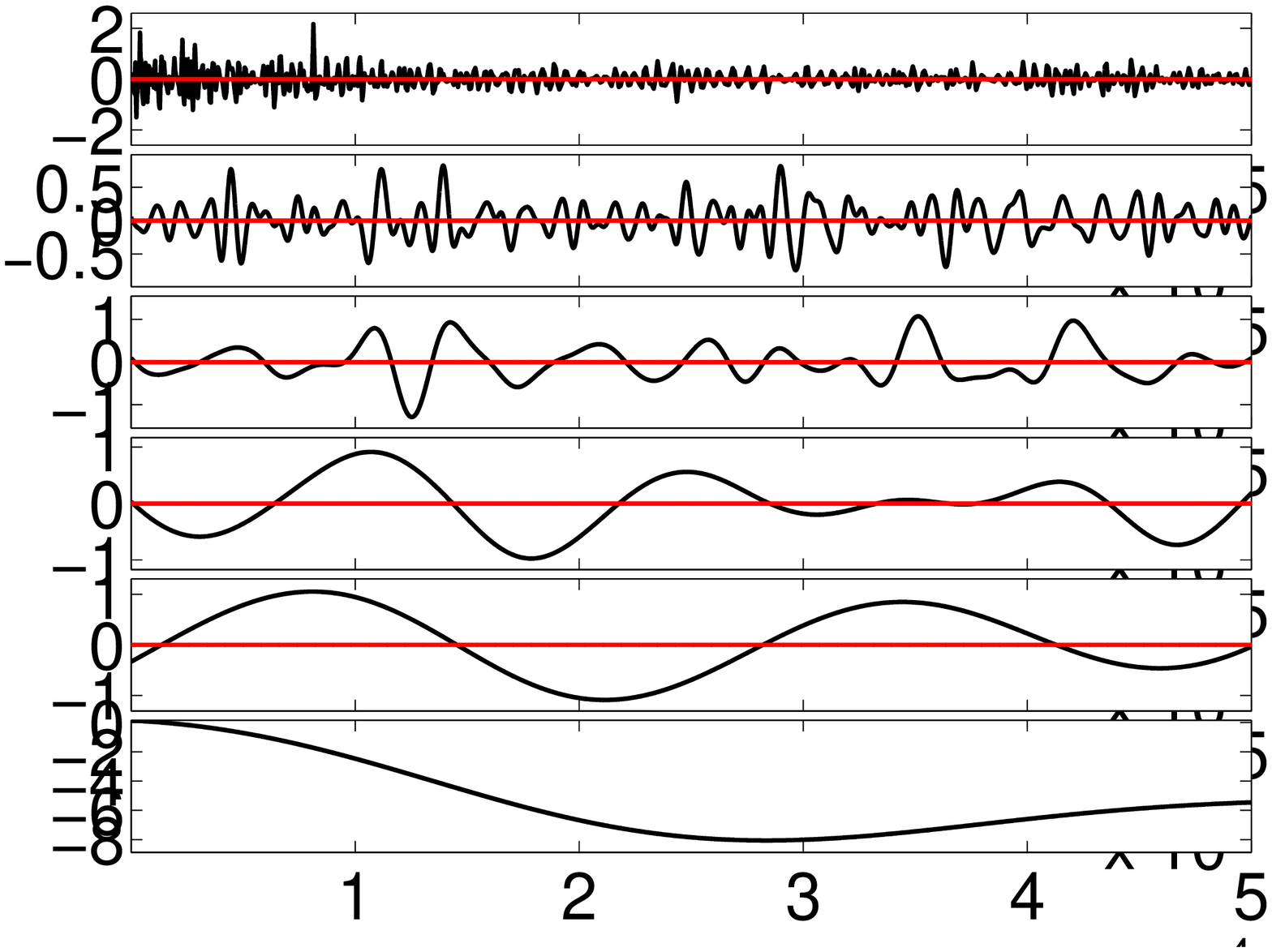}
                \caption{}
                \label{fig:Paleo_IF_Decomp}
        \end{subfigure}
        ~ 
        \begin{subfigure}[b]{0.5\textwidth}
                \centering
                \includegraphics[width=\textwidth]{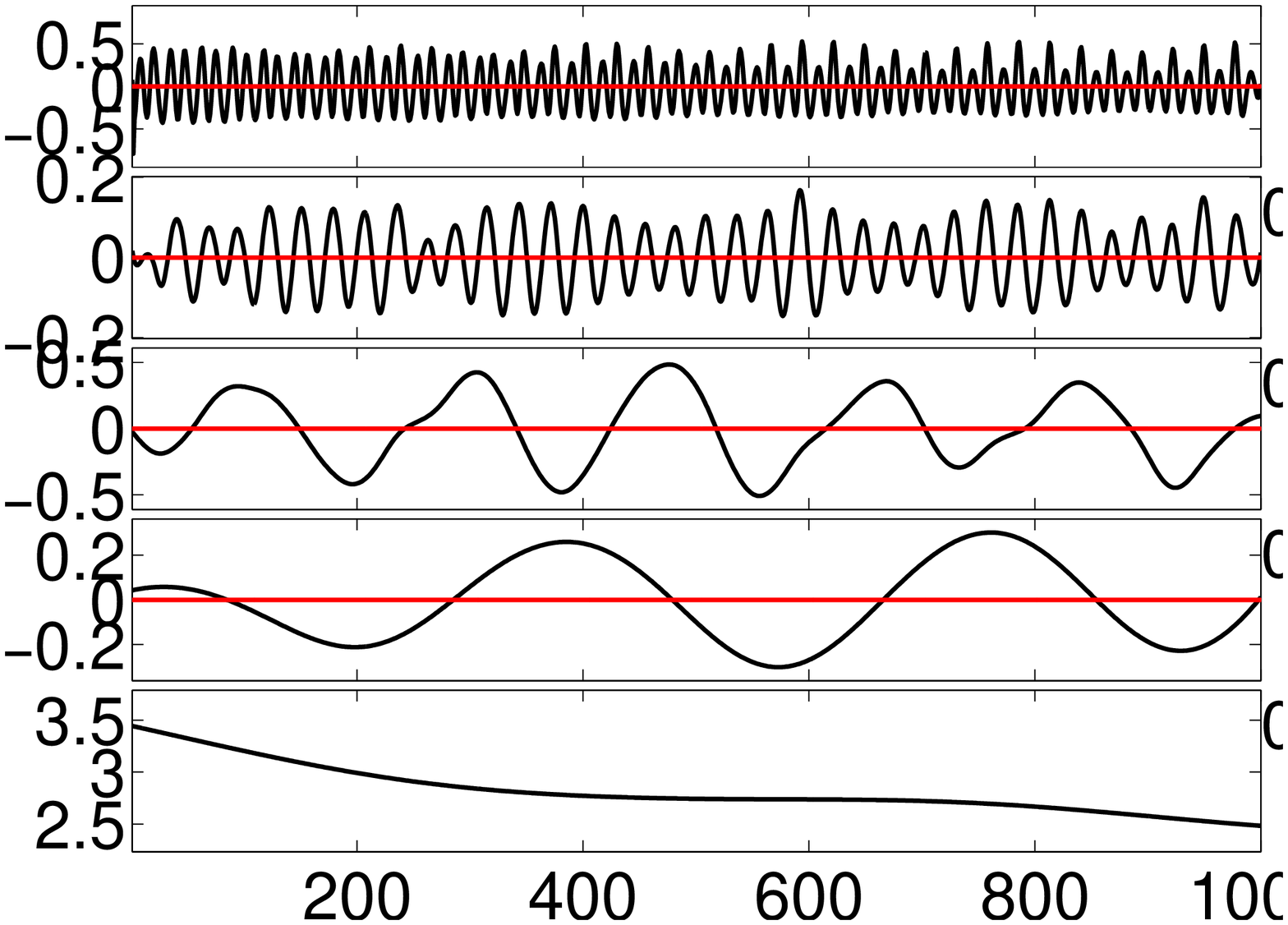}
                \caption{}
                \label{fig:LOD_IF_Decomp}
        \end{subfigure}
        \caption{(\subref{fig:Paleo_IF_Decomp}) Vostok temperature IF decomposition. (\subref{fig:LOD_IF_Decomp}) length--of--day IF decomposition.}\label{fig:IF_decomp}
\end{figure}

Regarding the Vostok temperature we set NIMFs to 100, for alpha we use the default value \verb"ave". If we use the default value of \verb"Xi" equal to 1.6 the method decomposes the low frequency component in several IMFs. Therefore we increase \verb"Xi" up to 3 using the script \verb"Settings_IF('IF.NIMFs',100,'IF.Xi',3);".
We obtain a decomposition which contains 5 IMFs and a remainder, Figure \ref{fig:Paleo_IF_Decomp}.

For the LOD data we use the \verb"Almost_min"  option for the parameter alpha and we set NIMFs to 100. Also in this example we have to increase \verb"Xi" up to 3. We use the script \verb"Settings_IF('IF.NIMFs',100,'IF.alpha','Almost_min','IF.Xi',3);". We obtain a decomposition which contains 4 IMFs and a remainder, Figure \ref{fig:LOD_IF_Decomp}.

\section{Conclusion and outlook}\label{sec:Outlook}

In this paper we provide the reader with a brief overview of methods for the decomposition of nonstationary signals. We address more in details the question of how to decompose such signals into simpler components using iterative methods which do not require any a priori assumptions.

From a research point of view, while the EEMD is still lacking a rigorous mathematical analysis, the IF algorithm has been extensively studied and analyzed in 1D and higher dimensions \cite{cicone2014adaptive,cicone2017multidimensional,cicone2017numerical,cicone2017spectral}. However it is important to point out that the Iterative Filtering method may fail to meaningfully decompose signals whose instantaneous frequencies vary consistently over time like, for instance, in a chirp. For this reason the Adaptive Local Iterative Filtering (ALIF) method has been proposed in \cite{cicone2014adaptive}. It generalizes the IF algorithm allowing for meaningful decompositions of any kind of nonstationary signals even with wide and quick changes in the instantaneous frequencies, like chirps. However ALIF mathematical understanding is far from being complete \cite{cicone2017spectral}. More research has to be done in this direction.

\section*{Acknowledgements}

The author's research was supported by Istituto Nazionale di Alta Matematica (INdAM) ``INdAM Fellowships in Mathematics and/or Applications cofunded by Marie Curie Actions'', PCOFUND-GA-2009-245492 INdAM-COFUND Marie Sklodowska Curie Integration Grants.

The author is deeply grateful to Haomin Zhou, a great researcher and wonderful person. He contributed substantially to this work and to the author career with many suggestions and pieces of advice he gave to the author over the years.

\end{document}